\documentclass[11pt]{article}
\usepackage{amssymb}
\usepackage{amsmath}
\usepackage[all]{xy}

\title{\bf Generalized Tilting Modules With Finite Injective
Dimension
\thanks{2000 {\it Mathematics Subject
Classification}. 16E10, 16E30.}
\thanks{{\it Key words and phrases}. generalized tilting modules, injective
dimension, $U$-limit dimension, submodule-closed, Wakamatsu Tilting
Conjecture.}}
\author{Zhaoyong Huang\thanks{\small \it E-mail address: huangzy@nju.edu.cn}\\
{\small \it Department of Mathematics, Nanjing University,}\\
{\small \it Nanjing 210093, People's Republic of China}\\}
\usepackage{amssymb}
\date{}
\begin{document}
\baselineskip=18pt \maketitle

\begin{abstract}
Let $R$ be a left noetherian ring, $S$ a right noetherian ring and
$_RU$ a generalized tilting module with $S={\rm End}(_RU)$. The
injective dimensions of $_RU$ and $U_S$ are identical provided both
of them are finite. Under the assumption that the injective
dimensions of $_RU$ and $U_S$ are finite, we describe when the
subcategory $\{ {\rm Ext}_S^n(N, U)|N$ is a finitely generated right
$S$-module$\}$ is submodule-closed. As a consequence, we obtain a
negative answer to a question posed by Auslander in 1969. Finally,
some partial answers to Wakamatsu Tilting Conjecture are given.
\end{abstract}

\vspace{1cm}

\centerline{\bf 1. Introduction}

\vspace{0.2cm}

Let $R$ be a ring. We use Mod $R$ (resp. Mod $R ^{op}$) to denote
the category of left (resp. right) $R$-modules, and use mod $R$
(resp. mod $R ^{op}$) to denote the category of finitely generated
left $R$-modules (resp. right $R$-modules).

We define gen$^*(_RR)=\{ X \in {\rm mod}\ R|$there exists an exact
sequence $\cdots \to P_i \to \cdots \to P_1 \to P_0 \to X \to 0$ in
mod $R$ with  $P_i$ projective for any $i \geq 0 \}$ (see [W2]). A
module $_RU$ in mod $R$ is called {\it selforthogonal} if ${\rm
Ext}_R^i(_RU, {_RU})=0$ for any $i \geq 1$.

\vspace{0.2cm}

{\bf Definition 1.1}$^{{\rm[W2]}}$ A selforthogonal module $_RU$ in
gen$^*(_RR)$ is called a {\it generalized tilting module} (sometimes
it is also called a {\it Wakamatsu tilting module}, see [BR]) if
there exists an exact sequence:
$$0 \to {_{R}R} \to U_0 \to U_1 \to \cdots \to U_i \to \cdots$$
such that: (1)  $U_i \in$add$_{R}U$ for any $i \geq 0$, where
add$_RU$ denotes the full subcategory of mod $R$ consisting of all
modules isomorphic to direct summands of finite sums of copies of
$_{R}U$, and (2) after applying the functor Hom$_{R}(\ , {_RU})$ the
sequence is still exact.

\vspace{0.2cm}

Let $R$ and $S$ be any rings. Recall that a bimodule $_RU_S$ is
called a {\it faithfully balanced bimodule} if the natural maps $R
\rightarrow {\rm End}(U_S)$ and $S \rightarrow {\rm End}(_RU)^{op}$
are isomorphisms. By [W2] Corollary 3.2, we have that $_{R}U _{S}$
is faithfully balanced and selforthogonal with $_RU \in$gen$^*(_RR)$
and $U_S \in$gen$^*(S_S)$ if and only if $_{R}U$ is generalized
tilting with $S =$End$(_{R}U)$ if and only if $U_{S}$ is generalized
tilting with $R =$End$(U_{S})$.

Let $R$ and $S$ be Artin algebras and $_RU$ a generalized tilting
module with $S={\rm End}(_{R}U)$. Wakamatsu proved in [W1] Theorem
that the projective (resp. injective) dimensions of $_RU$ and $U_S$
are identical provided both of them are finite. The result on the
projective dimensions also holds true when $R$ is a left noetherian
ring and $S$ is a right noetherian ring, by using an argument
similar to that in [W1]. In this case, $_RU_S$ is a tilting bimodule
of finite projective dimension ([M] Proposition 1.6). However,
because there is no duality available, Wakamatsu's argument in [W1]
does not work on the injective dimensions over noetherian rings. So,
it is natural to ask the following questions: When $R$ is a left
noetherian ring and $S$ is a right noetherian ring, (1) Do the
injective dimensions of $_RU$ and $U_S$ coincide provided both of
them are finite? (2) If one of the injective dimensions of $_RU$ and
$U_S$ is finite, is the other also finite?

The answer to the first question is positive if one of the following
conditions is satisfied: (1) $_RU_S={_RR_R}$ ([Z] Lemma A); (2) $R$
and $S$ are Artin algebras ([W1] Theorem); (3) $R$ and $S$ are
two-sided noetherian rings and $_RU$ is $n$-Gorenstein for all $n$
([H2] Proposition 17.2.6). In this paper, we show in Section 2 that
the answer to this question is always positive.

By the positive answer to the first question, the second question is
equivalent to the following question: Are the injective dimensions
of $_{R}U$ and $U_S$ identical? The above result means that the
answer to this question is positive provided that both dimensions
are finite. On the other hand, for Artin algebras, the positive
answer to the second question is equivalent to the validity of
Wakamtsu Tilting Conjecture ({\bf WTC}). This conjecture states that
every generalized tilting module with finite projective dimension is
tilting, or equivalently, every generalized tilting module with
finite injective dimension is cotilting. Moreover, {\bf WTC} implies
the validity of the Gorenstein Symmetry Conjecture ({\bf GSC}),
which states that the left and right self-injective dimensions of
$R$ are identical (see [BR]). In Section 4, we give some partial
answers to question (2). Let $R$ and $S$ be two-sided artinian rings
and $_RU$ a generalized tilting module with $S={\rm End}(_{R}U)$. We
prove that if the injective dimension of $U_S$ is equal to $n$ and
the $U$-limit dimension of each of the first $(n-1)$-st terms is
finite, then the injective dimension of $_RU$ is also equal to $n$.
Thus it trivial that the injective dimension of $U_S$ is at most 1
if and only if that of $_RU$ is at most 1. We remark that for an
Artin algebra $R$, it is well known that the right self-injective
dimension of $R$ is at most 1 if and only if the left self-injective
dimension of $R$ is at most 1 (see [AR3] p.121). In addition, we
prove that the left and right injective dimensions of $_RU$ and
$U_S$ are identical if $_RU$ (or $U_S$) is quasi Gorenstein, that
is, {\bf WTC} holds for quasi Gorenstein modules.

For an $(R-S)$-bimodule $_RU_S$ and a positive integer $n$, we
denote $\mathcal{E}_n(U_S)=\{ M\in {\rm mod}\ R|M={\rm Ext}_S^n(N,
U)$ for some $N \in$mod $S^{op} \}$. For a two-sided noetherian ring
$R$, Auslander showed in [A] Proposition 3.3 that any direct summand
of a module in $\mathcal{E}_1(R_R)$ is still in
$\mathcal{E}_1(R_R)$. He then asked whether any submodule of a
module in $\mathcal{E}_1(R_R)$ is still in $\mathcal{E}_1(R_R)$.
Recall that a full subcategory $\mathcal{X}$ of mod $R$ is said to
be {\it submodule-closed} if any non-zero submodule of a module in
$\mathcal{X}$ is also in $\mathcal{X}$. Then the above Auslander's
question is equivalent to the following question: Is
$\mathcal{E}_1(R_R)$ submodule-closed? In Section 3, under the
assumption that $R$ is a left noetherian ring, $S$ is a right
noetherian ring and $_RU$ is a generalized tilting module with
$S={\rm End}(_{R}U)$ and the injective dimensions of $_RU$ and $U_S$
being finite, we give some necessary and sufficient conditions for
$\mathcal{E}_n(U_S)$ being submodule-closed. As a consequence, we
construct some examples to illustrate that neither
$\mathcal{E}_1(R_R)$ nor $\mathcal{E}_2(R_R)$ are submodule-closed
in general, by which we answer the above Auslander's question
negatively.

Throughout this paper, $R$ is a left noetherian ring, $S$ is a right
noetherian ring (unless stated otherwise) and $_RU$ is a generalized
tilting module with $S={\rm End}(_RU)$. For a module $A$ in Mod $R$
(resp. Mod $S^{op}$), we use {\it l.}id$_R(A)$, {\it l.}fd$_R(A)$
and {\it l.}pd$_R(A)$ (resp. {\it r.}id$_S(A)$, {\it r.}fd$_S(A)$
and {\it r.}pd$_S(A)$) to denote the injective dimension, flat
dimension and projective dimension of $_RA$ (resp. $A_S$),
respectively.

\vspace{0.5cm}

\centerline{\bf 2. Some homological dimensions}

\vspace{0.2cm}

In this section, we study the relations among the $U$-limit
dimension (which was introduced in [H2]) of an injective module $E$,
the flat dimension of Hom$(U, E)$ and the injective dimension of
$U$. Then we show that {\it l.}id$_R(U)$={\it r.}id$_S(U)$ provided
both of them are finite.

The following result is [W2] Corollary 3.2.

\vspace{0.2cm}

{\bf Proposition 2.1} {\it The following statements are equivalent}.

(1) $_RU$ {\it is a generalized tilting module with} $S={\rm
End}(_RU)$.

(2) $U_S$ {\it is a generalized tilting module with} $R={\rm
End}(U_S)$.

(3) $_RU_S$ {\it is a faithfully balanced and selforthogonal
bimodule}.

\vspace{0.2cm}

We use add-lim$_RU$ to denote the subcategory of Mod $R$ consisting
of all modules isomorphic to direct summands of a direct limit of a
family of modules in which each is a finite direct sum of copies of
$_RU$ (see [H2]).

\vspace{0.2cm}

{\bf Proposition 2.2} (1) {\it Let} $V \in$add-lim$_RU$. {\it Then}
Ext$_R^i(U^{(I)}, V)=0$ {\it for any index set} $I$ {\it and} $i\geq
1$.

(2) Ext$_R^i(U^{(I)}, U^{(J)})=0$ {\it for any index sets} $I$, $J$
{\it and} $i\geq 1$.

\vspace{0.2cm}

{\it Proof.} (1) It is well known that for any $i \geq 1$,
Ext$_R^i(U^{(I)}, V)\cong$Ext$_R^i(U, V)^I$. Since $_RU$ is finitely
generated and selforthogonal and $V \in$add-lim$_RU$, it follows
easily from [S1] Theorem 3.2 that Ext$_R^i(U, V)=0$ and so
Ext$_R^i(U, V)^I=0$.

(2) Because a direct sum of a family of modules is a special kind of
a direct limit of these modules, (2) follows from (1) trivially.
\hfill{$\square$}

\vspace{0.2cm}

{\bf Definition 2.3}$^{{\rm [H2]}}$ For a module $A$ in Mod $R$, if
there exists an exact sequence $\cdots \to U_{n} \to \cdots \to
U_{1} \to U_{0} \to A \to 0$ in Mod $R$ with  $U_{i}\in$add-lim$_RU$
for any $i \geq 0$, then we define the $U$-{\it limit dimension} of
$A$, denoted by $U$-lim.dim$_R(A)$, as inf$\{ n|$there exists an
exact sequence $0 \to U_{n} \to \cdots \to U_{1} \to U_{0} \to A \to
0$ in Mod $R$ with  $U_{i}\in$add-lim$_RU$ for any $0 \leq i \leq n
\}$. We set $U$-lim.dim$_R(A)$ infinity if no such an integer
exists.

\vspace{0.2cm}

{\it Remark.} It is well known that a module over any ring is flat
if and only if it is direct limit of a family of finitely generated
free modules. So, putting $_RU={_RR}$, a module in Mod $R$ is flat
if and only if it is in add-lim$_{R}R$; in this case, the dimension
defined as in Definition 2.3 is just the flat dimension of modules.

\vspace{0.2cm}

For a module $A$ in Mod $R$ (resp. Mod $S^{op}$), we denote either
of Hom$_{R}(_{R}U _{S}, {_RA})$ and Hom$_{S}(_{R}U _{S}, {A_S})$ by
$^*A$.

\vspace{0.2cm}

{\bf Lemma 2.4} {\it Let} $_RE$ {\it be an injective} $R$-{\it
module. Then l.}fd$_S(^*E)=U$-lim.dim$_R(E)$.

\vspace{0.2cm}

{\it Proof.} The result was proved in [H2] Lemma 17.3.1 when $R$ and
$S$ are two-sided noetherian rings. The proof in [H2] remains valid
in the setting here, we omit it. \hfill{$\square$}

\vspace{0.2cm}

{\it Remark.} It is not difficult to see from the proof of [H2]
Lemma 17.3.1 that for any injective $R$-module $E$, there exists an
exact sequence:
$$\cdots \to U_i \to \cdots \to U_1 \to U_0 \to E \to 0$$
in Mod $R$ with  $U_i$ in add-lim$_RU$ for any $i \geq 0$. So
$U$-lim.dim$_R(E)$ (finite or infinite) always exists for any
injective $R$-module $E$.

\vspace{0.2cm}

{\bf Proposition 2.5} {\it Let} $_RE$ {\it be an injective} $R$-{\it
module. If l.}id$_R(U)=n<\infty$ {\it and} {\it
l.}fd$_S(^*E)<\infty$, {\it then} {\it l.}fd$_S(^*E)\leq n$.

\vspace{0.2cm}

{\it Proof.} Suppose {\it l.}fd$_S(^*E)=m<\infty$. Then there exists
an exact sequence:
$$0 \to F_m \to S^{(I_{m-1})} \to \cdots \to S^{(I_1)} \to S^{(I_0)}
\to {^*E} \to 0 \eqno{(1)}$$ in Mod $S$ with $F_m$ flat and  $I_i$
an index set for any $0 \leq i \leq m-1$. By [CE] Chapter VI,
Proposition 5.3, we have that
$${\rm Tor}_{j}^{S}(U , {^*E}) \cong {\rm
Hom}_{R}({\rm Ext}_{S}^{j}(U , U), E)=0$$ for any $j \geq 1$. So, by
applying the functor $U\otimes _S-$ to the exact sequence (1), we
get the following exact sequence:
$$0 \to U\otimes _SF_m \to U\otimes _SS^{(I_{m-1})} \to \cdots \to
U\otimes _SS^{(I_1)} \to U\otimes _SS^{(I_0)} \to U\otimes _S{^*E}
\to 0.$$ By Proposition 2.1 and [S2] p.47, we have that $U\otimes
_S{^*E} \cong {\rm Hom}_{R}({\rm Hom}_{S}(U, U), E) \cong E$. So we
get the following exact sequence:
$$0 \to K_m \buildrel {d_m} \over \longrightarrow
U^{(I_{m-1})} \buildrel {d_{m-1}} \over \longrightarrow \cdots
\buildrel {d_2} \over \longrightarrow U^{(I_1)} \buildrel {d_1}
\over \longrightarrow U^{(I_0)} \buildrel {d_0} \over
\longrightarrow E \to 0 \eqno{(2)}$$ where $K_m=U\otimes _SF_m$.
Because $F_m$ is a flat $S$-module, it is a direct limit of finitely
generated free $S$-modules. So $K_m=U\otimes _SF_m \in$add-lim$_RU$.
Hence, by Proposition 2.2(1), Ext$_R^j(U^{(I_i)}, K_m)=0$ for any
$j\geq 1$ and $0 \leq i \leq m-1$.

Since {\it l.}id$_R(U)=n$, {\it l.}id$_R(K_m)\leq n$ by [S1] Theorem
3.2. If $m>n$, then Ext$_R^m(E, K_m)=0$. It follows from the exact
sequence (2) that Ext$_R^1(K_{m-1}, K_m)=0$, where
$K_{m-1}$=Coker$d_m$. Thus the sequence $0 \to K_m \buildrel {d_m}
\over \longrightarrow U^{(I_{m-1})} \to K_{m-1} \to 0$ splits and
$U^{(I_{m-1})}\cong K_m \bigoplus K_{m-1}$. In addition, we get an
exact sequence:
$$0 \to K_{m-1} \to U^{(I_{m-2})}\buildrel {d_{m-2}} \over \longrightarrow \cdots
\buildrel {d_2} \over \longrightarrow U^{(I_1)} \buildrel {d_1}
\over \longrightarrow U^{(I_0)} \buildrel {d_0} \over
\longrightarrow E \to 0$$ with $K_{m-1}, U^{(I_{m-2})}, \cdots,
U^{(I_0)} \in$add-lim$_RU$. Then $U$-lim.dim$_R(E) \leq m-1$. But
{\it l.}fd$_S(^*E)=U$-lim.dim$_R(E)$ by Lemma 2.4. Consequently we
conclude that {\it l.}fd$_S(^*E) \leq${\it l.}id$_R(U)$.
\hfill{$\square$}

\vspace{0.2cm}

We also need the following result, which is [H2] Lemma 17.2.4.

\vspace{0.2cm}

{\bf Lemma 2.6} (1) {\it r.}id$_{S}(U)$=sup$\{${\it l.}fd$_{S}(^*E)|
_{R}E$ {\it is injective}$\}$. {\it Moreover, r.}id$_{S}(U)$={\it
l.}fd$_{S}(^*Q)$ {\it for any injective cogenerator} $_{R}Q$ {\it
for} Mod $R$.

(2) {\it l.}id$_{R}(U)$=sup$\{${\it r.}fd$_{R}(^* E')|E'_{S}$ {\it
is injective}$\}$. {\it Moreover, l.}id$_{R}(U)$={\it
r.}fd$_{R}(^*Q')$ {\it for any injective cogenerator} $Q'_{S}$ for
Mod $S^{op}$.

\vspace{0.2cm}

We are now in a position to prove one of the main results in this
paper.

\vspace{0.2cm}

{\bf Theorem 2.7} {\it l.}id$_R(U)$={\it r.}id$_S(U)$ {\it provided
both of them are finite}.

\vspace{0.2cm}

{\it Proof.} Let $_RQ$ be an injective cogenerator for Mod $R$.
Assume that {\it l.}id$_R(U)=n<\infty$ and {\it
r.}id$_S(U)=m<\infty$. Then {\it l.}fd$_{S}(^*Q)=m$ by Lemma 2.6. So
$m=${\it l.}fd$_{S}(^*Q)\leq${\it l.}id$_{R}(U)=n$ by Proposition
2.5. Dually, we may prove $n \leq m$. We are done. \hfill{$\square$}

\vspace{0.2cm}

{\bf Definition 2.8}$^{[{\rm AB}]}$ Let $\mathcal{X}$ be a full
subcategory of Mod $R$. For a module $A$ in Mod $R$, if there exists
an exact sequence $\cdots \to T_{n} \to \cdots \to T_{1} \to T_{0}
\to A \to 0$ in Mod $R$ with  $T_{i}\in\mathcal{X}$ for any $i \geq
0$, then we define $\mathcal{X}$-resol.dim$_R(A)$=inf$\{ n|$there
exists an exact sequence $0 \to T_{n} \to \cdots \to T_{1} \to T_{0}
\to A \to 0$ in Mod $R$ with  $T_{i}\in\mathcal{X}$ for any $0 \leq
i \leq n \}$. We set $\mathcal{X}$-resol.dim$_R(A)$ infinity if no
such an integer exists.

\vspace{0.2cm}

We use Add$_RU$ to denote the full subcategory of Mod $R$ consisting
of all modules isomorphic to direct summands of sums of copies of
$_RU$. Compare the following result with Lemma 2.4.

\vspace{0.2cm}

{\bf Lemma 2.9} {\it Let} $_RE$ {\it be an injective} $R$-{\it
module. Then l.}pd$_S(^*E)$=Add$_{R}U$-resol.dim$_{R}(E)$.

\vspace{0.2cm}

{\it Proof.} We first prove that Add$_{R}U$-resol.dim$_{R}(E)
\leq${\it l.}pd$_S(^*E)$. Without loss of generality, assume that
{\it l.}pd$_S(^*E)=m<\infty$. Then there exists an exact sequence:
$$0 \to Q_m \to Q_{m-1} \to \cdots \to Q_1 \to Q_0
\to {^*E} \to 0$$ in Mod $S$ with  $Q_i$ projective for any $0 \leq
i \leq m$. Then by using an argument similar to that in the proof of
Proposition 2.5, we get that Add$_{R}U$-resol.dim$_{R}(E) \leq m$.

We next prove that {\it
l.}pd$_S(^*E)\leq$Add$_{R}U$-resol.dim$_{R}(E)$. Assume that
Add$_{R}U$-resol.dim$_{R}(E)$ \linebreak $=m<\infty$. Then there
exists an exact sequence: $$0 \to U_{m} \to \cdots \to U_{1} \to
U_{0} \to E \to 0 \eqno{(3)}$$ in Mod $R$ with  $U_{i}\in$Add$_{R}U$
for any $0 \leq i \leq m$. By Proposition 2.2(2), we have that
$^*U_i \in$Add$_SS$ (that is, $^*U_i$ is a projective left
$S$-module) and Ext$_R^j(U, U_i)=0$ for any $j \geq 1$ and $0 \leq i
\leq m$. So by applying the functor Hom$_R({_RU}, -)$ to the exact
sequence (3), we get the following exact sequence:
$$0 \to {^*U_{m}} \to \cdots \to {^*U_{1}} \to
{^*U_{0}} \to {^*E} \to 0$$ in Mod $S$ with  $^*U_{i}$ left
$S$-projective for any $0 \leq i \leq m$, and hence {\it
l.}pd$_S(^*E)\leq m$. \hfill{$\square$}

\vspace{0.2cm}

{\it Remark.} (1) Put $_RU={_RR}$. Then
Add$_RU$-resol.dim$_R(A)$={\it l.}pd$_R(A)$ for any $A \in$Mod $R$.

(2) Because a direct sum of a family of modules is a special kind of
a direct limit of these modules, for any $A \in$Mod $R$, we have
that $U$-lim.dim$_R(A) \leq$Add$_RU$-resol.dim$_R(A)$ if both of
them exist.

(3) It is not difficult to see from the proof of Lemma 2.9 that for
any injective $R$-module $E$, there exists an exact sequence:
$$\cdots \to U_i \to \cdots \to U_1 \to U_0 \to E \to 0$$
in Mod $R$ with  $U_i$ in Add$_RU$ for any $i \geq 0$. So
Add$_{R}U$-resol.dim$_{R}(E)$ (finite or infinite) always exists for
any injective $R$-module $E$.

\vspace{0.2cm}

The proof of Proposition 2.5 in fact proves the following more
general result.

\vspace{0.2cm}

{\bf Proposition 2.10} {\it Let} $_RE$ {\it be an injective}
$R$-{\it module. If l.}id$_R(U)=n<\infty$ {\it and} {\it
l.}fd$_S(^*E)<\infty$, {\it then} Add$_{R}U$-resol.dim$_{R}(E)\leq
n$ ({\it equivalently, l.}pd$_S(^*E)\leq n)$.

\vspace{0.2cm}

{\bf Theorem 2.11} {\it Let} $_RE$ {\it be an injective} $R$-{\it
module. If l.}id$_R(U)=n<\infty$, {\it then the following statements
are equivalent}.

(1) {\it l.}fd$_S(^*E)<\infty$.

(2) {\it l.}pd$_S(^*E)<\infty$.

(3) $U$-lim.dim$_R(E)<\infty$.

(4) Add$_{R}U$-resol.dim$_{R}(E)<\infty$.

(5) {\it l.}fd$_S(^*E)\leq n$.

(6) {\it l.}pd$_S(^*E)\leq n$.

(7) $U$-lim.dim$_R(E)\leq n$.

(8) Add$_{R}U$-resol.dim$_{R}(E)\leq n$.

\vspace{0.2cm}

{\it Proof.} The implications that $(6)\Rightarrow (2) \Rightarrow
(1)$ and $(6)\Rightarrow (5)\Rightarrow (1)$ are trivial. The
implication of $(1)\Rightarrow (6)$ follows from Proposition 2.10.
By Lemma 2.4, we have $(1)\Leftrightarrow (3)$ and
$(5)\Leftrightarrow (7)$. By Lemma 2.9, we have $(2)\Leftrightarrow
(4)$ and $(6)\Leftrightarrow (8)$. \hfill{$\square$}

\vspace{0.2cm}

As an application of the obtained results, we get the following
corollary, which gives some equivalent conditions that {\it
l.}id$_R(U)=n$ implies {\it r.}id$_S(U)=n$.

\vspace{0.2cm}

{\bf Corollary 2.12} {\it Let} $_RQ$ {\it be an injective
cogenerator for} Mod $R$. {\it If l.}id$_R(U)=n(<\infty)$, {\it then
the following statements are equivalent.}

(1) {\it r.}id$_S(U)=n$.

(2) {\it One of l.}fd$_S(^*Q)$, {\it l.}pd$_S(^*Q)$,
$U$-lim.dim$_R(Q)$ {\it and} Add$_{R}U$-resol.dim$_{R}(Q)$ {\it is
finite}.

\vspace{0.2cm}

{\it Proof.} Let $_RQ$ be an injective cogenerator for Mod $R$. Then
by Lemmas 2.6(1), 2.4 and 2.9, we have that {\it r.}id$_S(U)$={\it
l.}fd$_S(^*Q)$=$U$-lim.dim$_R(Q)
\leq$Add$_{R}U$-resol.dim$_{R}(Q)$={\it l.}pd$_S(^*Q)$. Now the
equivalence of (1) and (2) follows easily from Theorems 2.11 and
2.7. \hfill{$\square$}

\vspace{0.5cm}

\newpage

\centerline{\bf 3. Submodule-closure of $\mathcal{E}_n(U_S)$}

\vspace{0.2cm}

In this section, we study Auslander's question mentioned in Section
1 in a more general situation.

\vspace{0.2cm}

{\bf Lemma 3.1} {\it For any injective module} $_RE$ {\it and any
non-negative integer} $t$, {\it l.}fd$_S(^*E) \leq t$ {\it if and
only if} Hom$_R({\rm Ext}_S^{t+1}(N, U), E)=0$ {\it for any module}
$N \in$mod $S^{op}$.

\vspace{0.2cm}

{\it Proof.} It is easy by [CE] Chapter VI, Proposition 5.3.
\hfill{$\square$}

\vspace{0.2cm}

For a module $A \in$ mod $R$ and a non-negative integer $n$, we say
that the {\it grade} of $A$ with respect to $_{R}U$, written as
grade$_UA$, is at least $n$ if Ext$_{R}^i(A, U)=0$ for any $0 \leq i
< n$. We say that the {\it strong grade} of $A$ with respect to
$_RU$, written as s.grade$_UA$, is at least $n$ if grade$_UB \geq n$
for all submodules $B$ of $A$ (see [H2]). Assume that
$$0 \to {_RU} \to E_0 \buildrel {\alpha _0} \over \longrightarrow
E_1 \buildrel {\alpha _1} \over \longrightarrow \cdots \to E_i
\buildrel {\alpha _i} \over \longrightarrow \cdots$$ is a minimal
injective resolution of $_RU$.

\vspace{0.2cm}

{\bf Lemma 3.2} Let $n$ be a positive integer and $m$ an integer
with $m \geq -n$. Then the following statements are equivalent.

(1) $U$-lim.dim$_R(\bigoplus _{i=0}^{n-1}E_{i})\leq n+m$.

(2) s.grade$_U{\rm Ext}_S^{n+m+1}(N, U) \geq n$ {\it for any} $N
\in$mod $S^{op}$.

(3) {\it l.}fd$_S(^*E_i) \leq n+m$ {\it for any} $0 \leq i \leq
n-1$.

\vspace{0.2cm}

{\it Proof.} This conclusion has been proved in [H2] Lemma 17.3.2
when $R$ and $S$ are two-sided noetherian rings. The argument there
remains valid in the setting here, we omit it. \hfill{$\square$}

\vspace{0.2cm}

For a module $A$ in mod $R$ (resp. mod $S^{op}$), we call
Hom$_{R}(_{R}A, {_{R}U _{S}})$ (resp. Hom$_{S}(A_{S}, {_{R}U_{S}})$)
the dual module of $A$ with respect to $_{R}U_{S}$, and denote
either of these modules by $A^*$. For a homomorphism $f$ between
$R$-modules (resp. $S^{op}$-modules), we put $f^*={\rm Hom}(f,
{_{R}U_{S}})$. We use $\sigma _{A}: A \rightarrow A^{**}$ via
$\sigma _{A}(x)(f)=f(x)$ for any $x \in A$ and $f \in A^*$ to denote
the canonical evaluation homomorphism. $A$ is called $U$-{\it
torsionless} (resp. $U$-{\it reflexive}) if $\sigma _{A}$ is a
monomorphism (resp. an isomorphism).

\vspace{0.2cm}

{\bf Definition 3.3} ([H3]) Let $\mathcal{X}$ be a full subcategory
of mod $R$. $\mathcal{X}$ is said to have the $U$-{\it torsionless
property} (resp. the $U$-{\it reflexive property}) if each module in
$\mathcal{X}$ is $U$-torsionless (resp. $U$-reflexive).

\vspace{0.2cm}

We denote $^{\bot}_RU=\{ M\in {\rm mod}\ R|{\rm Ext}_R^i(M,
{_RU})=0$ for any $i \geq 1\}$ and $^{\bot _n}_RU=\{ M\in {\rm mod}\
R|{\rm Ext}_R^i(M, {_RU})=0$ for any $1 \leq i \leq n\}$ (where $n$
is a positive integer). A module $M$ in mod $R$ is said to have {\it
generalized Gorenstein dimension zero} (with respect to ${_RU_S}$),
denoted by G-dim$_U(M)=0$, if the following conditions are
satisfied: (1) $M$ is $U$-reflexive, and (2) $M\in {^{\bot}_RU}$ and
$M^*\in {^{\bot}U_S}$. Symmetrically, we may define the notion of a
module in mod $S^{op}$ having generalized Gorenstein dimension zero
(with respect to ${_RU_S}$) (see [AR2]). We use $\mathcal{G}_U$ to
denote the full subcategory of mod $R$ consisting of the modules
with generalized Gorenstein dimension zero. It is trivial that
$^{\bot}_RU\supseteq \mathcal{G}_U$.

\vspace{0.2cm}

{\bf Proposition 3.4} ([H3] Proposition 2.3) {\it The following
statements are equivalent}.

(1) $_R^{\bot}U$ {\it has the} $U$-{\it torsionless property}.

(2) $_R^{\bot}U$ {\it has the} $U$-{\it reflexive property}.

(3) $^{\bot}_RU=\mathcal{G}_U$.

\vspace{0.2cm}

For any $n\geq 0$, we denote $\mathcal{H}_n(_RU)=\{ M\in {\rm mod}\
R| {\rm Ext}_R^i(M, {_RU})=0$ for any $i\geq 0$ with $i\neq n\}$
([W2]]).

\vspace{0.2cm}

{\bf Lemma 3.5} {\it If} $_R^{\bot}U$ {\it has the} $U$-{\it
torsionless property, then} $\mathcal{H}_n(_RU)\subseteq
\mathcal{E}_n(U_S)$ {\it for any} $n \geq 1$.

\vspace{0.2cm}

{\it Proof.} It follows from Proposition 3.4 and [H4] Lemma 3.3.
\hfill{$\square$}

\vspace{0.2cm}

{\bf Lemma 3.6} {\it Assume that} $_R^{\bot _n}U$ {\it has the}
$U$-{\it torsionless property, where} $n$ {\it is a positive
integer. If} $A$ {\it is a non-zero module in} mod $R$ {\it with}
grade$_UA \geq n$, {\it then} grade$_UA=n$.

\vspace{0.2cm}

{\it Proof.} Let $0\neq A \in$mod $R$ with grade$_UA \geq n$. If
grade$_UA>n$, then $A^*=0$ and $A \in {_R^{\bot _n}U}$. Since
$_R^{\bot _n}U$ has the $U$-torsionless property, $A$ is
$U$-torsionless and $A \hookrightarrow A^{**}=0$, which is a
contradiction. Thus grade$_UA=n$. \hfill{$\square$}

\vspace{0.2cm}

For any $n\geq 0$, we denote $\overline{\mathcal{H}_n(_RU)}=\{ M\in
{\rm mod}\ R|$ any non-zero submodule of $M$ is in
$\mathcal{H}_n(_RU)\}$. It is clear that
$\overline{\mathcal{H}_n(_RU)}\subseteq \mathcal{H}_n(_RU)$. We are
now in a position to give the main result in this section.

\vspace{0.2cm}

{\bf Theorem 3.7} {\it If l.}id$_R(U) \leq n$ {\it and} $_R^{\bot
_n}U$ {\it has the} $U$-{\it torsionless property, where} $n$ {\it
is a positive integer, then the following statements are
equivalent.}

(1) $U$-lim.dim$_{R}(\bigoplus _{i=0}^{n-1}E_{i})\leq n-1$.

(2) $\mathcal{E}_n(U_S)$ {\it is submodule-closed and}
$\mathcal{E}_n(U_S)=\mathcal{H}_n(_RU)$.

(3) $\mathcal{E}_n(U_S)$ {\it is submodule-closed and}
$\mathcal{E}_n(U_S)=\overline{\mathcal{H}_n(_RU)}$.

{\it Proof.} Since $_R^{\bot _n}U$ has the $U$-torsionless property,
$\overline{\mathcal{H}_n(_RU)}\subseteq \mathcal{H}_n(_RU) \subseteq
\mathcal{E}_n(U_S)$ by Lemma 3.5. So the implication of $(3)
\Rightarrow (2)$ is trivial.

$(1) \Rightarrow (3)$ Assume that $U$-lim.dim$_{R}(\bigoplus
_{i=0}^{n-1}E_{i})\leq n-1$ and $M$ is any non-zero module in
$\mathcal{E}_n(U_S)$. Then s.grade$_UM \geq n$ by Lemma 3.2.

Let $A$ be any non-zero submodule of $M$ in mod $R$. Then grade$_UA
\geq n$. By Lemma 3.6, grade$_UA=n$. In addition, {\it l.}id$_R(U)
\leq n$, so $A \in \mathcal{H}_n(_RU)$ and $M \in
\overline{\mathcal{H}_n(_RU)}$. It follows that
$\mathcal{E}_n(U_S)\subseteq \overline{\mathcal{H}_n(_RU)}$ and
$\mathcal{E}_n(U_S)=\overline{\mathcal{H}_n(_RU)}$.

Notice that $\overline{\mathcal{H}_n(_RU)}$ is clearly
submodule-closed, so $\mathcal{E}_n(U_S)$ is also submodule-closed.

$(2) \Rightarrow (1)$ We first prove $U$-lim.dim$_{R}(E_0)\leq n-1$.
If $U$-lim.dim$_{R}(E_0)>n-1$, then {\it l.}fd$_{S}(^*E_0)
>n-1$ by Lemma 2.4. So by Lemma 3.1, there exists a module $N
\in$mod $S^{op}$ such that Hom$_R({\rm Ext}_S^n(N, U), E_0)\neq 0$.
Hence there exists a non-zero homomorphism $f: {\rm Ext}_S^n(N, U)
\to E_0$. Since $_RU$ is essential in $E_0$, $f^{-1}(_RU)$ is a
non-zero submodule of ${\rm Ext}_S^n(N, U)$. By assumption,
$\mathcal{H}_n(_RU)=\mathcal{E}_n(U_S)$ and $\mathcal{E}_n(U_S)$ is
submodule-closed. So $f^{-1}(_RU) \in
\mathcal{E}_n(U_S)(=\mathcal{H}_n(_RU))$ and hence
$[f^{-1}(_RU)]^*=0$, which is a contradiction. Consequently,
$U$-lim.dim$_{R}(E_0)\leq n-1$.

We next prove $U$-lim.dim$_{R}(E_1)\leq n-1$ (note: at this moment,
$n \geq 2$). If $U$-lim.dim$_{R}(E_1)$ \linebreak $>n-1$, then {\it
l.}fd$_{S}(^*E_1)>n-1$ by Lemma 2.4. So by Lemma 3.1, there exists a
module $N_1 \in$mod $S^{op}$ such that Hom$_R({\rm Ext}_S^n(N_1, U),
E_1)\neq 0$. Hence there exists a non-zero homomorphism $f_1: {\rm
Ext}_S^n(N_1, U) \to E_1$. Since ${\rm Ker}\alpha _1$ is essential
in $E_1$, $f_1^{-1}({\rm Ker}\alpha _1)$ is a non-zero submodule of
${\rm Ext}_S^n(N_1, U)$. By assumption, $f_1^{-1}({\rm Ker}\alpha
_1) \in \mathcal{E}_n(U_S)$. Since {\it
l.}fd$_{S}(^*E_0)=U$-lim.dim$_{R}(E_0)\leq n-1$,
Hom$_R(f_1^{-1}({\rm Ker}\alpha _1), E_0)=0$ by Lemma 3.1.

From the exact sequence $0 \to {_RU} \to E_0 \to {\rm Ker}\alpha _1
\to 0$ we get the following exact sequence:
$${\rm Hom}_R(f_1^{-1}({\rm Ker}\alpha _1), E_0) \to
{\rm Hom}_R(f_1^{-1}({\rm Ker}\alpha _1), {\rm Ker}\alpha _1) \to
{\rm Ext}_R^1(f_1^{-1}({\rm Ker}\alpha _1), {_RU}).$$ Since
$f_1^{-1}({\rm Ker}\alpha _1) \in
\mathcal{E}_n(U_S)(=\mathcal{H}_n(_RU))$, ${\rm
Ext}_R^1(f_1^{-1}({\rm Ker}\alpha _1), {_RU})=0$. So ${\rm
Hom}_R(f_1^{-1}({\rm Ker}\alpha _1),$ \linebreak ${\rm Ker}\alpha
_1)=0$, which is a contradiction. Hence we conclude that
$U$-lim.dim$_{R}(E_1)\leq n-1$.

Continuing this process, we get that $U$-lim.dim$_{R}(E_{i})\leq
n-1$ for any $0 \leq i \leq n-1$. \hfill{$\square$}

\vspace{0.2cm}

If {\it r.}id$_S(U) \leq n$, then $_R^{\bot _n}U$ has the
$U$-reflexive property by [HT] Theorem 2.2. So by Theorem 3.7, we
have the following

\vspace{0.2cm}

{\bf Corollary 3.8} {\it If l.}id$_R(U)$={\it r.}id$_S(U) \leq n$,
{\it then the following statements are equivalent.}

(1) $U$-lim.dim$_{R}(\bigoplus _{i=0}^{n-1}E_{i})\leq n-1$.

(2) $\mathcal{E}_n(U_S)$ {\it is submodule-closed and}
$\mathcal{E}_n(U_S)=\mathcal{H}_n(_RU)$.

(3) $\mathcal{E}_n(U_S)$ {\it is submodule-closed and}
$\mathcal{E}_n(U_S)=\overline{\mathcal{H}_n(_RU)}$.

\vspace{0.2cm}

Let $R$ be a two-sided noetherian ring. Recall that $R$ is called an
{\it Iwanaga-Gorenstein ring} if the injective dimensions of $_RR$
and $R_R$ are finite. Also recall that $R$ is said to satisfy the
{\it Auslander condition} if the flat dimension of the $(i+1)$-st
term in a minimal injective resolution of $_RR$ is at most $i$ for
any $i \geq 0$, and $R$ is called {\it Auslander-Gorenstein} if it
is Iwanaga-Gorenstein and satisfies the Auslander condition (see
[Bj]). It is well known that any commutative Iwanaga-Gorenstein ring
is Auslander-Gorenstein.

The following corollary gives a positive answer to the Auslander's
question for Auslander-Gorenstein rings and so in particular for
commutative Iwanaga-Gorenstein rings.

\vspace{0.2cm}

{\bf Corollary 3.9} {\it If} $R$ {\it is an Auslander-Gorenstein
ring with self-injective dimension} $n$, {\it then}
$\mathcal{E}_n(R_R)$ {\it is submodule-closed}.

\vspace{0.2cm}

{\it Proof.} Notice that $R$-lim.dim$_R(A)$={\it l.}fd$_R(A)$ for
any $A \in$Mod $R$, so our assertion follows from Corollary 3.8.
\hfill{$\square$}

\vspace{0.2cm}

Assume that $X \in$mod $S^{op}$ and there exists an exact sequence
$H_1 \buildrel {g} \over \longrightarrow H_0 \to X \to 0$ in mod
$S^{op}$. We denote $A=$Coker$g^*$. The following result is a
generalization of [HT] Lemma 2.1. The proof here is similar to that
in [HT], we omit it.

\vspace{0.2cm}

{\bf Lemma 3.10} {\it Let} $X$, $A$, $H_0$ {\it and} $H_1$ {\it be
as above. Assume that} $H_0$ {\it and} $H_1$ {\it are} $U$-{\it
reflexive}.

(1) {\it If} $H_i^* \in {^{\bot _{i+1}}_RU}$ {\it for} $i=0, 1$,
{\it then we have the following exact sequence}:
$$0 \to {\rm Ext}_R^1(A, U) \to X \buildrel {\sigma _X}
\over \longrightarrow X^{**} \to {\rm Ext}_R^2(A, U) \to 0.$$

(2) {\it If} $H_i \in {^{\bot _{2-i}}U_S}$ {\it for} $i=0, 1$, {\it
then we have the following exact sequence}:
$$0 \to {\rm Ext}_S^1(X, U) \to A \buildrel {\sigma _A}
\over \longrightarrow A^{**} \to {\rm Ext}_S^2(X, U) \to 0.$$

\vspace{0.2cm}

{\bf Lemma 3.11} {\it Let} $\mathcal{X}$ {\it be a full subcategory
of} $^{\bot}U_S$ {\it which has the} $U$-{\it reflexive property
and} $X$ {\it a module in} mod $S^{op}$. {\it If}
$\mathcal{X}$-resol.dim$_{S}(X)=n\ (\geq 1)$, {\it then}
grade$_U{\rm Ext}_S^n(X, U)\geq 1$; {\it if furthermore} $n \geq 2$
{\it and} $Y^* \in {^{\bot _2}_RU}$ {\it for any} $Y \in
\mathcal{X}$, {\it then} grade$_U{\rm Ext}_S^n(X, U)\geq 2$.

\vspace{0.2cm}

{\it Proof.} Assume that $\mathcal{X}$-resol.dim$_{S}(X)=n\ (\geq
1)$. Then there exists an exact sequence: $$0 \to X_{n} \buildrel
{d_n} \over \longrightarrow \cdots \to X_{1} \to X_{0} \to X \to 0$$
in mod $S^{op}$ with  $X_i \in \mathcal{X}\ (\subseteq
{^{\bot}U_S})$ for any $0 \leq i \leq n$. Set $N={\rm Coker}d_n$.
Then ${\rm Ext}_S^1(N, U) \cong {\rm Ext}_S^n(X, U)$.

Consider the following commutative diagram with exact rows:
$$\xymatrix{& 0 \ar[r] & X_n \ar[r]^{d_n}
\ar[d]^{\sigma _{X_n}} & X_{n-1} \ar[r]
\ar[d]^{\sigma _{X_{n-1}}} & N \ar[r] & 0\\
0 \ar[r] & [{\rm Ext}_S^n(X, U)]^* \ar[r] & X_n^{**}
\ar[r]^{d_n^{**}} & X_{n-1}^{**} & &  }$$ Because $\mathcal{X}$ has
the $U$-reflexive property, both $\sigma _{X_n}$ and $\sigma
_{X_{n-1}}$ are isomorphisms. So we have that $[{\rm Ext}_S^n(X,
U)]^*=0$ and grade$_U{\rm Ext}_S^n(X, U)\geq 1$.

If $n \geq 2$ and $Y^* \in {^{\bot _2}_RU}$ for any $Y \in
\mathcal{X}$, by applying Lemma 3.10(1) to the exact sequence $0 \to
X_n \buildrel {d_n} \over \longrightarrow X_{n-1} \to N \to 0$, we
then get the following exact sequence:
$$0 \to {\rm Ext}_R^1({\rm Ext}_S^n(X, U), U) \to N \buildrel {\sigma _N}
\over \longrightarrow N^{**} \to {\rm Ext}_R^2({\rm Ext}_S^n(X, U),
U) \to 0.$$ Because $X_{n-2} \in \mathcal{X}$ and $\mathcal{X}$ has
the $U$-reflexive property, $X_{n-2}$ is $U$-reflexive. Then $N$ is
$U$-torsionless for it is isomorphic to a submodule of $X_{n-2}$. So
$\sigma _N$ is monic and ${\rm Ext}_R^1({\rm Ext}_S^n(X, U),$
\linebreak $U)=0$. Hence we conclude that grade$_U{\rm Ext}_S^n(X,
U)\geq 2$. \hfill{$\square$}

\vspace{0.2cm}

For a non-negative integer $t$, a module $N$ in mod $S^{op}$ is said
to have {\it generalized Gorenstein dimension} at most $t$ (with
respect to $_RU_S$), denoted by G-dim$_U(N) \leq t$, if there exists
an exact sequence $0 \to N_t \to \cdots \to N_1 \to N_0 \to N \to 0$
in mod $S^{op}$ with G-dim$_U(N_i)=0$ for any $0 \leq i \leq t$ (see
[AR2]).

\vspace{0.2cm}

{\bf Lemma 3.12} {\it Let} $n \leq 2$. {\it If l.}id$_R(U)$={\it
r.}id$_S(U) \leq n$, {\it then}
$\mathcal{E}_n(U_S)=\mathcal{H}_n(_RU)$.

\vspace{0.2cm}

{\it Proof.} Because {\it l.}id$_R(U)$={\it r.}id$_S(U) \leq n$,
both $_R^{\bot _n}U$ and $^{\bot _n}U_S$ have the $U$-reflexive
property by [HT] Theorem 2.2. It follows from Lemma 3.5 that
$\mathcal{H}_n(_RU)\subseteq \mathcal{E}_n(U_S)$.

Assume that $n=1$ and $0\neq M \in \mathcal{E}_1(U_S)$. Let $E'_0$
be the injective envelope of $U_S$. Because {\it
r.}fd$_R(^*E'_0)\leq${\it l.}id$_R(U)\leq 1$ by assumption and Lemma
2.6, it follows from the symmetric statements of [H2] Theorem 17.5.5
that grade$_UM\geq 1$. By Lemma 3.6, grade$_UM=1$. Thus $M \in
\mathcal{H}_1(_RU)$ and $\mathcal{E}_1(U_S)\subseteq
\mathcal{H}_1(_RU)$. The case $n=1$ follows.

Assume that $n=2$ and $0\neq M \in \mathcal{E}_2(U_S)$. Then there
exists a module $N \in$mod $S^{op}$ such that $M={\rm Ext}_S^2(N,
U)(\neq 0)$. Because {\it l.}id$_R(U)$={\it r.}id$_S(U) \leq 2$, by
[HT] Theorem 3.5 we have that G-dim$_U(N) \leq 2$. If G-dim$_U(N)
<2$, then there exists an exact sequence $0 \to N_1 \to N_0 \to N
\to 0$ in mod $S^{op}$ with G-dim$_U(N_1)$=G-dim$_U(N_0)=0$. So
${\rm Ext}_S^2(N, U) \cong {\rm Ext}_S^1(N_1, U)=0$, which is a
contradiction. Hence we conclude that G-dim$_U(N)=2$. Then by Lemma
3.11, grade$_UM$=grade$_U{\rm Ext}_S^2(N, U) \geq 2$. By Lemma 3.6,
grade$_UM=2$. Thus $M \in \mathcal{H}_2(_RU)$ and
$\mathcal{E}_2(U_S)\subseteq \mathcal{H}_2(_RU)$. The case $n=2$
follows. \hfill{$\square$}

\vspace{0.2cm}

{\bf Proposition 3.13} {\it Let} $n \leq 2$. {\it Assume that} {\it
l.}id$_R(U)$={\it r.}id$_S(U) \leq n$.

(1) {\it If} $n=1$, {\it then} $\mathcal{E}_1(U_S)$ {\it is
submodule-closed if and only if} $E_0 \in$add-lim$_RU$ ({\it that
is,} $U$-lim.dim$_{R}(E_{0})=0$) {\it if and only if}
$U$-lim.dim$_{R}(E_{i}) \leq i$ {\it for} $i=0, 1$.

(2) {\it If} $n=2$, {\it then} $\mathcal{E}_2(U_S)$ {\it is
submodule-closed if and only if} $U$-lim.dim$_{R}(E_0 \bigoplus
E_{1})\leq 1$.

\vspace{0.2cm}

{\it Proof.} The former equivalence in (1) and the equivalence in
(2) follow from Lemma 3.12 and Theorem 3.7. Notice that
$U$-lim.dim$_{R}(E_{1})$={\it l.}fd$_S(^*E_1) \leq${\it
r.}id$_{S}(U)$ by Lemmas 2.4 and 2.6, then the latter equivalence in
(1) follows. \hfill{$\square$}

\vspace{0.2cm}

Let
$$0 \to {_RR} \to I_0 \to
I_1 \to \cdots \to I_i \to \cdots$$ be a minimal injective
resolution of $_RR$. Putting $_RU_S={_RR_R}$, by Proposition 3.13 we
immediately have the following

\vspace{0.2cm}

{\bf Corollary 3.14} {\it Let} $n \leq 2$. {\it Assume that
l.}id$_R(R)$={\it r.}id$_R(R) \leq n$.

(1) {\it If} $n=1$, {\it then} $\mathcal{E}_1(R_R)$ {\it is
submodule-closed if and only if} $I_0$ {\it is flat if and only if}
$R$ {\it is Auslander-Gorenstein}.

(2) {\it If} $n=2$, {\it then} $\mathcal{E}_2(R_R)$ {\it is
submodule-closed if and only if l.}fd$_{R}(I_0 \bigoplus I_{1})\leq
1$.

\vspace{0.2cm}

In the following, we give some examples to illustrate that neither
$\mathcal{E}_1(R_R)$ nor $\mathcal{E}_2(R_R)$ are submodule-closed
in general.

\vspace{0.2cm}

{\bf Example 3.15} Let $K$ be a field and $R$ a finite dimensional
$K$-algebra which is given by the quiver:
$$2 \longleftarrow 1 \longrightarrow 3.$$
Then $R$ is Iwanaga-Gorenstein with {\it l.}id$_R(R)$={\it
r.}id$_R(R)=1$ and {\it l.}fd$_R(I_0)=1$. By Corollary 3.14,
$\mathcal{E}_1(R_R)$ is not submodule-closed.

\vspace{0.2cm}

{\bf Example 3.16} Let $K$ be a field and $\Delta$ the quiver:
$$\xymatrix{1 \ar[r]^{\gamma} \ar[d]^{\alpha} & 3 \ar[d]^{\delta}\\
2 \ar[r]^{\beta} & 4 }$$ If $R=K\Delta /(\beta \alpha)$, then $R$ is
Iwanaga-Gorenstein with {\it l.}id$_R(R)$={\it r.}id$_R(R)=2$ and
{\it l.}fd$_R(E(P_4))$ \linebreak $=2$, where $E(P_4)$ is the
injective envelope of the indecomposable projective module
corresponding to the vertex 4. Since $P_4$ is a direct summand of
$_RR$, {\it l.}fd$_R(I_0)=2$. By Corollary 3.14,
$\mathcal{E}_2(R_R)$ is not submodule-closed.

\vspace{0.2cm}

It is clear that mod $R \supseteq \mathcal{E}_1(R_R) \supseteq
\mathcal{E}_2(R_R) \supseteq \cdots \supseteq \mathcal{E}_i(R_R)
\supseteq \cdots$. From the above argument we know that
$\mathcal{E}_n(R_R)$ is submodule-closed for an Auslander-Gorenstein
ring $R$ with self-injective dimension $n$ for any $n \geq 1$, and
neither $\mathcal{E}_1(R_R)$ nor $\mathcal{E}_2(R_R)$ are
submodule-closed in general. However, we don't know whether
$\mathcal{E}_n(R_R)$ (where $n\geq 3$) is submodule-closed or not in
general.

\vspace{0.5cm}

\centerline{\bf 4. Wakamatsu tilting conjecture and (quasi)
Gorenstein modules}

\vspace{0.2cm}

Let $R$ be an Artin algebra. Recall that a module $_RT$ in mod $R$
is called a {\it tilting module} of finite projective dimension if
the following conditions are satisfied: (1) {\it
l}.pd$_R(T)<\infty$; (2) $_RT$ is selforthogonal; and (3) there
exists an exact sequence $0 \to R \to T_0 \to T_1 \to \cdots \to T_t
\to 0$ in mod $R$ with
 $T_i \in$add$_RT$ for any $0 \leq i \leq t$. The notion of
{\it cotilting modules} of finite injective dimension may be defined
dually. A generalized tilting module is not necessarily tilting or
cotilting. The following conjecture is called Wakamatsu Tilting
Conjecture ({\bf WTC}): Every generalized tilting module with finite
projective dimension is tilting, or equivalently, every generalized
tilting module with finite injective dimension is cotilting (see
[BR]). For Artin algebras $R$ and $S$ and a generalized tilting
module $_RU$ with $S$=End$(_RU)$, by Theorem 2.7 and the dual
results of [M] Theorem 1.5 and Proposition 1.6, we easily get the
following equivalent statements:

(1) {\bf WTC} holds.

(2) If one of {\it l.}id$_R(U)$ and {\it r.}id$_S(U)$ is finite,
then the other is also finite.

(3) {\it l.}id$_R(U)$={\it r.}id$_S(U)$.

The Gorenstein Symmetry Conjecture ({\bf GSC}) states that the left
and right self-injective dimensions of $R$ are identical for an
Artin algebra $R$ (see [BR]). It is trivial from the above
equivalent conditions that {\bf WTC} $\Rightarrow$ {\bf GSC}. As an
application of the results obtained in Section 2, we now give some
sufficient conditions for the validity of statement (2). In other
words, we establish some cases in which {\bf WTC} holds true.

\vspace{0.2cm}

{\bf Theorem 4.1} ([H1] Theorem) {\it Let} $R$ {\it and} $S$ {\it be
two-sided artinian rings and} $n$ {\it and} $m$ {\it positive
integers. If r.}id$_S(U) \leq n$ {\it and} grade$_U$Ext$_R^m(M,
U)\geq n-1$ {\it for any} $M \in$mod $R$, {\it then l.}id$_R(U) \leq
n+m-1$.

\vspace{0.2cm}

The following corollary is an immediate consequence of Theorems 4.1
and 2.7.

\vspace{0.2cm}

{\bf Corollary 4.2} {\it Let} $R$ {\it and} $S$ {\it be two-sided
artinian rings. Then r.}id$_S(U) \leq 1$ {\it if and only if} {\it
l.}id$_R(U)\leq 1$.

\vspace{0.2cm}

Let $$0 \to {_RU} \to E_0 \to E_1 \to \cdots \to E_i \to \cdots$$
and
$$0 \to {U_S} \to E'_0 \to E'_1 \to \cdots \to E'_i \to
\cdots$$ be minimal injective resolutions of $_RU$ and $U_S$,
respectively. The following Propositions 4..3 and 4.6 generalize
some results in [H1] and [AR4].

\vspace{0.2cm}

{\bf Proposition 4.3} {\it Let} $R$ {\it and} $S$ {\it be two-sided
artinian rings and} $n$ {\it a positive integer. If r.}id$_S(U)=n$
{\it and} $U$-lim.dim$_S(\bigoplus _{i=0}^{n-2}E'_i)<\infty$, then
{\it l.}id$_R(U)=n$.

\vspace{0.2cm}

{\it Proof.} Assume that $U$-lim.dim$_S(\bigoplus
_{i=0}^{n-2}E'_i)=r (<\infty)$. It follows from [H2] Lemma 17.3.2
that s.grade$_U$Ext$_R^{r+1}(M, U) \geq n-1$ for any $M \in$mod $R$.
By Theorem 4.1, {\it l.}id$_R(U) \leq r+n (<\infty)$. Thus {\it
l.}id$_R(U)$={\it r.}id$_S(U)=n$ by Theorem 2.7. \hfill{$\square$}

\vspace{0.2cm}

The following two results are cited from [H2].

\vspace{0.2cm}

{\bf Theorem 4.4} ([H2] Theorem 17.1.11) {\it Let} $R$ {\it and} $S$
{\it be two-sided noetherian rings. Then, for a positive integer}
$n$, {\it the following statements are equivalent}.

(1) s.grade$_U$Ext$_{R}^i(M, U) \geq i$ {\it for any} $M \in$mod $R$
{\it and} $1 \leq i \leq n$.

(1)$^{op}$ s.grade$_U$Ext$_{S}^i(N, U) \geq i$ {\it for any} $N
\in$mod $S^{op}$ {\it and} $1 \leq i \leq n$.

$_RU$ ({\it symmetrically} $U_S$) {\it is called an} $n$-{\it
Gorenstein module if one of the above equivalent conditions is
satisfied, and} $_RU$ ({\it symmetrically} $U_S$) {\it is called a
Gorenstein module if it is} $n$-{\it Gorenstein for all} $n$.

\vspace{0.2cm}

It follows from [H2] Corollary 17.1.12 that a two-sided noetherian
ring $R$ satisfies the Auslander condition if and only if $_RR$ is a
Gorenstein module.

\vspace{0.2cm}

{\bf Theorem 4.5} ([H2] Theorem 17.5.4) {\it Let} $R$ {\it and} $S$
{\it be two-sided noetherian rings. Then, for a positive integer}
$n$, {\it the following statements are equivalent}.

(1) s.grade$_U$Ext$_{S}^{i+1}(N, U) \geq i$ {\it for any} $N \in$mod
$S^{op}$ {\it and} $1 \leq i \leq n$.

(2) grade$_U$Ext$_{R}^i(M, U) \geq i$ {\it for any} $M \in$mod $R$
{\it and} $1 \leq i \leq n$.

$_RU$ {\it is called a quasi} $n$-{\it Gorenstein module if one of
the above equivalent conditions is satisfied, and} $_RU$ {\it is
called a quasi Gorenstein module if it is quasi} $n$-{\it Gorenstein
for all} $n$.

\vspace{0.2cm}

An ($n$-)Gorenstein module is clearly quasi ($n$-)Gorenstein. But
the conserve doesn't hold in general because the notion of
($n$)-Gorenstein modules is left-right symmetric by Theorem 4.4, and
that of quasi ($n$)-modules is not left-right symmetric even in the
case $_RU_S={_RR_R}$ (see [H2] Example 17.5.2).

\vspace{0.2cm}

{\bf Proposition 4.6} {\it Let} $R$ {\it and} $S$ {\it be two-sided
artinian rings. Then l.}id$_R(U)$={\it r.}id$_S(U)$ {\it provided
that} $_RU$ ({\it or} $U_S$) {\it is quasi Gorenstein}.

\vspace{0.2cm}

{\it Proof.} Let $_RU$ be a quasi Gorenstein module. By Theorem 4.5,
for any $i\geq 1$, we have that grade$_U$Ext$_{R}^i(M, U) \geq i$
for any $M \in$mod $R$ and s.grade$_U$Ext$_{S}^{i+1}(N, U) \geq i$
for any $N \in$mod $S^{op}$. Then it is easy to see from Theorem 4.1
that {\it l.}id$_R(U)<\infty$ if and only if {\it
r.}id$_S(U)<\infty$. Thus {\it l.}id$_R(U)$={\it r.}id$_S(U)$ by
Theorem 2.7. \hfill{$\square$}

\vspace{0.2cm}

Note that Proposition 4.6 generalizes [AR4] Corollary 5.5(b) which
asserts that {\it l.}id$_R(R)$=
\linebreak
{\it r.}id$_R(R)$ if $R$
is an Artin algebra satisfying the Auslander condition.

\vspace{0.2cm}

{\bf Conjecture 4.7} Let $R$ and $S$ be Artin algebras and $_RU$ a
generalized tilting module with $S={\rm End}(_RU)$. If $_RU$ is
(quasi) Gorenstein, then {\it l.}id$_R(U)$={\it r.}id$_S(U)<\infty$
(In fact, under our assumption it has been proved in Proposition 4.6
that {\it l.}id$_R(U)$={\it r.}id$_S(U)$).

\vspace{0.2cm}

Auslander and Reiten in [AR4] raised the following conjecture, which
we call Auslander Gorenstein Conjecture ({\bf AGC}): An Artin
algebra is Iwanaga-Gorenstein if it satisfies the Auslander
condition (in other words, an Artin algebra $R$ satisfies {\it
l.}id$_R(R)$={\it r.}id$_R(R)<\infty$ provided $_RR$ is a Gorenstein
module). It is trivial that this conjecture is situated between {\bf
Conjecture 4.7} and the famous Nakayama Conjecture ({\bf NC}), which
states that an Artin algebra $R$ is self-injective if each term in a
minimal injective resolution of $_RR$ is projective. That is, we
have the following implications: {\bf Conjecture 4.7} $\Rightarrow$
{\bf AGC} $\Rightarrow$ {\bf NC}.

Recall moreover the Generalized Nakayama Conjecture ({\bf GNC}):
Every indecomposable injective $R$-module occurs as the direct
summand of some term in a minimal injective resolution of $_RR$ for
an Artin algebra $R$. An equivalent version of {\bf GNC} is: For an
Artin algebra $R$ and every simple module $T \in$mod $R$, there
exists a non-negative integer $k$ such that Ext$_R^k(T, R) \neq 0$
(see [AR1]). It is well known that {\bf GNC} implies {\bf AGC}. We
now show the corresponding result for {\bf Conjecture 4.7}.

\vspace{0.2cm}

{\bf Proposition 4.8} {\it Let} $R$ {\it and} $S$ {\it be Artin
algebras and} $_RU$ {\it a generalized tilting module with} $S={\rm
End}(_RU)$. {\it If the following condition is satisfied: for every
simple module} $T \in$mod $R$, {\it there exists a negative-integer}
$k$ {\it such that} Ext$_R^k(T, {_RU})\neq 0$, {\it then} {\bf
Conjecture 4.7} {\it holds} for $R$.

\vspace{0.2cm}

{\it Proof.} Let $\{T_1$, $\cdots$, $T_t\}$ be the set of all
non-isomorphic simple modules in mod $R$. By assumption, for each
$T_i$ ($1 \leq i \leq t$), there exists a non-negative integer $k_i$
such that Ext$_R^{k_i}(T_i, {_RU})\neq 0$. It is easy to verify that
Hom$_R(T, E_j)\cong {\rm Ext}_R^j(T, {_RU})$ for any simple
$R$-module $T$ and $j \geq 0$. So Hom$_R(T_i, E_{k_i})\neq 0$ for
any $1 \leq i \leq t$ and hence  $E(T_i)$ (the injective envelope of
$T_i$) is isomorphic to a direct summand of $E_{k_i}$ for any $1
\leq i \leq t$.

Now suppose $_RU$ is quasi Gorenstein. Then by Theorem 4.5 and Lemma
3.2, we have that {\it l.}fd$_{S}(^*E_i) \leq i+1$ for any $i \geq
0$. So {\it l.}fd$_{S}(^*[E(T_i)]) \leq${\it
l.}fd$_{S}(^*E_{k_i})\leq k_i+1$ for any $1 \leq i \leq t$. Put
$E=\bigoplus _{i=1}^tE(T_i)$ and $k={\rm max}\{k_1, \cdots , k_t\}$.
Then $E$ is an injective cogenerator for Mod $R$ and {\it
l.}fd$_{S}(^*E) \leq k+1$. It follows from Lemma 2.6(1) that {\it
r.}id$_{S}(U)\leq k+1$. We are done. \hfill{$\square$}

\vspace{0.2cm}

Let $N$ be a module in mod $S^{op}$. Recall that an injective
resolution:
$$0 \to N \buildrel {\delta _0} \over \longrightarrow V_{0}
\buildrel {\delta _1} \over \longrightarrow V_{1} \buildrel {\delta
_2} \over \longrightarrow \cdots \buildrel {\delta _i} \over
\longrightarrow V_{i} \buildrel {\delta _{i+1}} \over
\longrightarrow \cdots$$ is called {\it ultimately closed} if there
exists a positive integer $n$ such that Im$\delta _n=\bigoplus
_{j=0}^mW_j$, where each $W_j$ is a direct summand of Im$\delta
_{i_j}$ with $i_j <n$. By [HT] Theorem 2.4, if $U_S$ has a
ultimately closed injective resolution (especially, if {\it
r.}id$_S(U)<\infty$), then $^{\bot}_RU$ has the $U$-reflexive
property and the condition in Proposition 4.8 is satisfied.

\vspace{0.5cm}

{\bf Acknowledgements} Part of the paper was written while the
author was staying at Universit${\rm \ddot{a}}$t Bielefeld supported
by the SFB 701 ``Spectral Structures and Topological Methods in
Mathematics". The author is grateful to Prof. Claus M. Ringel for
his warm hospitality and useful comments on this paper. The research
of the author was partially supported by Specialized Research Fund
for the Doctoral Program of Higher Education (Grant No. 20030284033,
20060284002) and NSF of Jiangsu Province of China (Grant No.
BK2005207). The author is grateful to the referee for the careful
reading and the valuable and detailed suggestions in shaping this
paper into its present version.

\end{document}